\documentclass[titlepage,twoside,12pt]{article}
\usepackage{amssymb}
\usepackage{amsfonts}
\begin{document}

{\LARGE \bf Where Infinitesimals Come From ...} \\ \\

{\bf Elem\'{e}r E ~Rosinger} \\ \\
{\small \it Department of Mathematics \\ and Applied Mathematics} \\
{\small \it University of Pretoria} \\
{\small \it Pretoria} \\
{\small \it 0002 South Africa} \\
{\small \it eerosinger@hotmail.com} \\ \\

\hspace*{6cm} {\it Dedicated to Marie-Louise Nykamp} \\ \\

{\bf Abstract} \\ 

The presence of infinitesimals is traced back to some of the most general algebraic structures, namely, semigroups, and in 
fact, magmas, [1], in which none of the structures of linear order, field, or the Archimedean property need to be present. 
Such a clarification of the basic structures from where infinitesimals can in fact emerge may prove to have a special
importance in Physics, as seen in [4-16]. The relevance of the deeper and simpler roots of infinitesimals, as they are 
given in Definitions 3.1 and 3.2, is shown by the close connection in Theorem 4.1 and Corollary 4.1 between the presence 
of infinitesimals and the non-Archimedean property, in the particular case of linearly ordered monoids, a case which, 
however, has a wide applicative interest. \\ \\

{\bf 1. Preliminaries} \\

Abraham Robinson, [3], considered it to be one of the important aspects of Nonstandard Analysis the first time rigorous
and comprehensive formulation of a theory of {\it infinitesimals}. As mentioned in the instructive historical survey at
the end of [3], Leibniz appears to be the first to introduce the idea of infinitesimals, and show their usefulness in
Calculus. And for the next two centuries, until in the second part of the 1800s Weierstrass introduced modern rigour into
the subject, Calculus had much been based on a variety of intuitive, rather than rigorous uses of infinitesimals popping
up in all kind of places and under any number of forms. In fact, in engineering or physics courses of Calculus, such loose
appeal to infinitesimals has gone on until more recently. \\

As it happens, at various times, a number of facts have not been understood quite clearly related to the status of
infinitesimals. \\

One such fact is that the Archimedean structure of the field $\mathbb{R}$ of usual real numbers does not allow the
presence of infinitesimals. This is the reason why after the reform introduced by Weierstrass there has no longer been a
place in Calculus for infinitesimals. \\
In this regard, Robinson's field $^*\mathbb{R}$ of nonstandard reals happens to be non-Archimedean, and as such, proves to
be able to accommodate infinitesimals. \\

However, in pursuing Nonstandard Analysis, Robinson had further goals in addition to obtaining a rigorous foundation for
infinitesimals. Indeed, among such goals was that $^*\mathbb{R}$ is a {\it linearly ordered field extension} of
$\mathbb{R}$. Furthermore, it was aimed that a good deal of the usual properties of $\mathbb{R}$ would automatically
remain valid for $^*\mathbb{R}$ as well, under the so called {\it transfer principle}. \\

A consequence of the above has been the tacit association of infinitesimals with :

\begin{itemize}

\item linear orders,

\item fields,

\item the Archimedean property.

\end{itemize}

As shown in this paper, however, the presence of infinitesimals can be traced back to far more general algebraic
structures, namely, semigroups, and in fact, magmas, [1], in which none of the above three structures need to be 
present. \\

Such a clarification of the basic structures from where infinitesimals can in fact emerge may prove to have a special 
importance in Physics, as seen in [4-16]. \\ 

The relevance of the deeper and simpler roots of infinitesimals, as they are given in Definitions 3.1 and 3.2, is shown by 
the close connection in Theorem 4.1 and Corollary 4.1 between the presence of infinitesimals and the non-Archimedean 
property in the particular case of linearly ordered monoids, a case which, however, has a wide applicative interest. \\ \\ 

{\bf 2. Rich and Complex Structure of the Set of \\
        \hspace*{0.45cm} Additive Subgroups} \\

It is useful to start by recalling the seldom considered and surprisingly rich and complex structure of the set of {\it additive subgroups} in
$^*\mathbb{R}$. \\

First, we recall that $\mathbb{R}$ is an additive subgroup of $^*\mathbb{R}$. Further, we can distinguish the following four subgroups in $^*\mathbb{R}$, namely \\

(2.1)~~~ $ \{ 0 \} ~\subsetneqq~ Mon ( 0 ) ~\subsetneqq~ Fin ( 0 ) ~\subsetneqq~ ^*\mathbb{R} $ \\

where $Mon ( 0 )$ denotes the set of infinitesimals, that is, the so called {\it monad} at $0 \in\, ^*\mathbb{R}$, while
$Fin ( 0 )$ denotes the set of finite elements $x \in\, ^*\mathbb{R}$, thus $Fin ( 0 ) = \mathbb{R} \bigoplus Mon ( 0 )$. \\

Clearly, when seen from $\mathbb{R}$, the subgroups (2.1) collapse to only two trivial instances, namely, $\{ 0 \}$ and $\mathbb{R}$ itself. \\

Now the important fact to note is that there are many more additive subgroups in $^*\mathbb{R}$, than listed in (2.1).
Indeed, let $\epsilon \in Mon ( 0 ),~ \epsilon > 0$, that is, a positive infinitesimal. Then, associated with this infinitesimal $\epsilon$ we obtain the
following infinitely many pair-wise disjoint additive subgroups in $^*\mathbb{R}$, namely \\

(2.2)~~~ $ \mathbb{R} \epsilon,~~ \mathbb{R} \epsilon^2,~~ \mathbb{R} \epsilon^n, \ldots ~\subsetneqq~ Mon ( 0 ) $ \\

And in fact, there are uncountably many of that type of pair-wise disjoint additive subgroups associated with the given infinitesimal $\epsilon$,
namely \\

(2.3)~~~ $ \mathbb{R} \epsilon^a ~\subsetneqq Mon ( 0 ),~~~ a \in \mathbb{R},~ a \geq 1 $ \\

Similarly, if we take any $X \in\, ^*\mathbb{R} \setminus Fin ( 0 ),~ X >0$, then associated with this infinitely large $X$ we obtain the following
infinitely many pair-wise disjoint additive subgroups in $^*\mathbb{R}$, namely \\

(2.4)~~~ $ \mathbb{R} X^a \subsetneqq\, ^*\mathbb{R},~~~ a \in \mathbb{R},~ a \geq 1 $ \\

Needless to say, the additive subgroups in $^*\mathbb{R}$ are far from being exhausted by those in (2.1) - (2.4). \\

As for the complexity of the relationships between various additive subgroups in $^*\mathbb{R}$, we can note the following. \\

Let  $\epsilon, \eta \in Mon ( 0 ),~ \epsilon, \eta > 0$, then \\

(2.5)~~~ $ \mathbb{R} \epsilon ~\bigcap~ \mathbb{R} \eta \neq \phi ~~~\Longleftrightarrow~~~
                          \mathbb{R} \epsilon = \mathbb{R} \eta ~~~\Longleftrightarrow~~~ \epsilon / \eta \in \mathbb{R} $ \\

and as is well known, the relation in the right hand of (2.5) is highly atypical among infinitesimals $\epsilon, \eta \in Mon ( 0 )$. \\

Similarly, let $X, Y \in\, ^*\mathbb{R} \setminus Fin ( 0 ),~ X, Y > 0$, then \\

(2.6)~~~ $ \mathbb{R} X ~\bigcap~ \mathbb{R} Y \neq \phi ~~~\Longleftrightarrow~~~
                          \mathbb{R} X = \mathbb{R} Y ~~~\Longleftrightarrow~~~ X / Y \in \mathbb{R} $ \\

where again, the relation in the right hand of (2.6) is highly atypical among infinitely large $X, Y \in\, ^*\mathbb{R} \setminus Fin ( 0 )$. \\

Let us now compare the above with the situation of additive subgroups of $\mathbb{R}$. \\

Given $x \in \mathbb{R},~ x > 0$, then \\

(2.7)~~~ $ \mathbb{R} x = \mathbb{R} $ \\

hence no trace of the rich complexity of additive subgroups such as in (2.1) - (2.7). \\

Let us consider another example with infinitesimals, one that is closely connected with the reduced power algebras, and in particular, with $^*\mathbb{R}$, see [4-16], namely the algebra $\mathbb{R}^\mathbb{N}$. \\

First we recall that we have the group isomorphism \\

(2.8)~~~ $ \mathbb{R} \ni x \longmapsto u_x = ( x, x, x, \ldots ) \in
                              {\cal U}_{\mathbb{R}^\mathbb{N}} \subset \mathbb{R}^\mathbb{N} $ \\

Further, in the algebra $\mathbb{R}^\mathbb{N}$ one can distinguish the following additive semigroups \\

(2.9)~~~ $ \{ 0 \} ~\subsetneqq~ {\cal I}_{\mathbb{R}^\mathbb{N}} ~\subsetneqq~ {\cal A}_{\mathbb{R}^\mathbb{N}}
                          ~\subsetneqq~ {\cal B}_{\mathbb{R}^\mathbb{N}} ~\subsetneqq~ \mathbb{R}^\mathbb{N} $ \\

where ${\cal I}_{\mathbb{R}^\mathbb{N}},~ {\cal A}_{\mathbb{R}^\mathbb{N}}$ and ${\cal B}_{\mathbb{R}^\mathbb{N}}$ are,
respectively, the set of sequences $x = ( x_0, x_1, x_2, \ldots ) \in \mathbb{R}^\mathbb{N}$, which converge to $0 \in
\mathbb{R}$, converge to some element in $\mathbb{R}$, respectively, are bounded. \\

In (2.9), in view of [4-16], one can see ${\cal I}_{\mathbb{R}^\mathbb{N}}$ as the monad of $0 \in \mathbb{R}^\mathbb{N}$,
that is, the set of infinitesimals in $\mathbb{R}^\mathbb{N}$, while $\mathbb{R}^\mathbb{N} \setminus {\cal
B}_{\mathbb{R}^\mathbb{N}}$ can be seen as the set of infinitely large elements in $\mathbb{R}^\mathbb{N}$. \\

Clearly, and unlike with $^*\mathbb{R}$, in the algebra $\mathbb{R}^\mathbb{N}$ it is {\it not} the case that \\

$~~~~~~ \epsilon ~\mbox{infinitesimal},~ \epsilon \neq 0 ~~~\Longrightarrow~~~ 1 / \epsilon ~\mbox{infinitely large} $ \\

since it may happen that $1 / \epsilon$ is not even defined. Similarly, it is {\it not} the case that \\

$~~~~~~ X ~\mbox{infinitely large} ~~~\Longrightarrow~~~ 1 / X ~\mbox{infinitesimal} $ \\

since it may happen that $1 / X$ is not even defined. \\

Similar with the situation in (2.2) - (2.6), and in fact, with an increased richness and complexity, one obtains the 
structure of the set of additive subgroups in $\mathbb{R}^\mathbb{N}$. \\ \\  

{\bf 3. Defining Infinitesimals} \\

The sharp contrast seen in section 2 between the richness and complexity of additive semigroups in $^*\mathbb{R}$ and
$\mathbb{R}^\mathbb{N}$, and on the other hand, in $\mathbb{R}$, suggests the following definition \\

{\bf Definition 3.1} \\

A semigroup $( E, * )$ is called a {\it hyperspace}, if and only if it contains a sub-semigroup $F \subsetneqq E$,
together with an infinite set ${\cal I}$ of pair-wise disjoint sub-semigroups $I \subsetneqq F$. \\

In such a case the elements of the set \\

(3.1)~~~ $ I_{\cal I} \,=\, \bigcup_{I \in {\cal I}}\, I $ \\

are called {\it infinitesimals}, when the set ${\cal I}$ is maximal with its respective property.

\hfill $\Box$ \\

Here, and in the sequel, two sub-semigroups are called disjoint, if and only if their intersection is void, or it is a set of 
idempotent elements. \\

Clearly, the above definition can be extended to {\it magmas}, [1], namely \\

{\bf Definition 3.2} \\

A magma $( E, * )$ is called a {\it hyperspace}, if and only if it contains a sub-magma $F \subsetneqq E$, together with
an infinite set ${\cal I}$ of pair-wise disjoint sub-magmas $I \subsetneqq F$. \\

In such a case the elements of the set \\

(3.2)~~~ $ I_{\cal I} \,=\, \bigcup_{I \in {\cal I}}\, I $ \\

are called {\it infinitesimals}, when the set ${\cal I}$ is maximal with its respective property. 

\hfill $\Box$ \\

Similar with above, here, and in the sequel, two sub-magmas are called disjoint, if and only if their intersection is void, 
or it is a set of idempotent elements. \\

{\bf Remark 3.1} \\

The above two definitions do not make use of any partial, let alone, linear order. Equally, they do not make use of but of one single algebraic
operation, unlike in the case of algebras or fields, where at least two operations, namely, addition and multiplication are involved. Finally, they do
not make use of any kind of Archimedean or non-Archimedean property. \\

As for the algebraic operation $*$ involved in the above two definitions, it takes the place of addition, rather than
multiplication, as suggested by the example in section 2. \\ \\

{\bf 4 Examples} \\

Let us start with two familiar and somewhat different versions of the concept of non-Archimedean structure within the
setting of partially ordered monoids. Namely, let $( E, +, \leq )$ be a partially ordered monoid, thus we have satisfied \\

(4.1)~~~ $ x, y \in E_+ ~~\Longrightarrow~~ x + y \in E_+ $ \\

where $E_+ = \{ x \in E ~|~ x \geq 0 \}$. \\

A first intuitive version of the Archimedean condition, suggested in case $\leq$ is a {\it linear} order on $E$, is \\

(4.2)~~~ $ \exists~~ u \in E_+ ~:~ \forall~~ x \in E ~:~
                            \exists~~ n \in \mathbb{N} ~:~ n u \geq x $ \\

Here however is an alternative condition used in the literature when $\leq$ may be a partial order on $E$ \\

(4.3)~~~ $ \begin{array}{l}
                \forall~~ x \in E_+ ~:~ \\
                 ~~~~ x = 0 ~\Longleftrightarrow~
                          \left ( \begin{array}{l}
                                       \exists~~ y \in E_+ ~: \\
                                       \forall~~ n \in \mathbb{N} ~: \\
                                        ~~~~ n x \leq y
                                   \end{array} \right )
          \end{array} $ \\ \\

where clearly the implication "$\Longrightarrow$" is trivial, and thus condition (4.3) is equivalent with \\

(4.4)~~~ $ \forall~~ x \in E_+ ~:~
               \mathbb{N} x ~~\mbox{is bounded above}~~
                  \Longrightarrow ~~ x = 0 $ \\

{\bf Lemma 4.1} \\

We have the implication (4.2) ~~$\Longrightarrow$~~ (4.3). \\

{\bf Proof} \\

Assume that (4.3), hence (4.4) does not hold, then \\

$~~~~~~ \exists~~ x \in E_+ ~:~
               \mathbb{N} x ~~\mbox{is bounded above, and}~~ x \neq 0 $ \\

thus \\

$~~~~~~ \exists~~ u \in E_+,~ x \in E ~:~ \forall~~ n \in \mathbb{N} ~:~ n u \leq x $ \\

and (4.1) is contradicted.

\hfill $\Box$ \\

As for the converse implication (4.3) ~~$\Longrightarrow$~~ (4.2), we have \\

{\bf Lemma 4.2} \\

If $( E, +, \leq )$ is a {\it linearly} ordered monoid, then (4.3) ~~$\Longrightarrow$~~ (4.2). \\

{\bf Proof} \\

Assume indeed that (4.2) does not hold, then \\

$~~~~~~ \forall~~ x \in E_+ ~:~ \exists~~ y \in E ~:~
                            \forall~~ n \in \mathbb{N} ~:~ n x \ngeqslant y $ \\

and since $\leq$ is a linear order on $E$, we have \\

$~~~~~~ \forall~~ x \in E_+ ~:~ \exists~~ y \in E ~:~
                            \forall~~ n \in \mathbb{N} ~:~ n x \leq y $ \\

Obviously, we can assume that $y \in E_+$, thus (4.4) is contradicted. \\

{\bf Theorem 4.1} \\

A nontrivial linearly ordered monoid $( E, +, \leq )$ which is non-Archimedean in the sense of (4.2), is a hyperspace. \\

{\bf Proof.} \\

Assume that contrary to (4.2), we have \\

(4.5)~~~ $ \forall~~ u \in E_+ ~:~ \exists~~ x \in E ~:~
                            \forall~~ n \in \mathbb{N} ~:~ n u \leq x $ \\

and we note that $x \geq u$, thus in particular $x \in E_+$. \\

Let us take $u_1 > 0$. Then (4.5) gives $x_1 \geq u_1 > 0$, such that \\

(4.6)~~~ $ I_1 = \mathbb{Z} u_1 \leq x_1 $ \\

We take now $u_2 > x_1$, and as above, we obtain $x_2 \geq u_2$, such that \\

(4.7)~~~ $ I_2 = \mathbb{Z} u_2 \leq x_2 $ \\

Continuing the procedure, we obtain \\

(4.8)~~~ $ 0 < u_1 \leq x_1 < u_2 \leq x_2 < \ldots $ \\

(4.9)~~~ $ I_1 = \mathbb{Z} u_1 \leq x_1,~~ I_2 = \mathbb{Z} u_2 \leq x_2, \ldots $ \\

We show now that \\

(4.10)~~~ $ I_1 \cap I_2 = \{ 0 \} $ \\

Indeed, let $y \in I_1 \cap I_2,~$, then \\

(4.11)~~~ $ y = n_1 u_1 = n_2 u_2,~ y \neq 0$ \\

for some $n_1, n_2 \in \mathbb{Z}$. But (4.8), (4.9) give \\

(4.12)~~~ $ y = n_1 u_1 \leq x_1 < u_2 $ \\

hence $n_2 u_2 = y < u_2$ which means \\ 

(4.13)~~~ $ n_2 < 0 $ \\ 

and thus (4.11) yields \\ 

(4.14)~~~ $ n_1 < 0$ \\ 

It follows that \\ 

()~~~ $ -y = n_1 u_1 = n_2 u_2 \in I_1 \cap I_2 $ \\

hence, from the start, we can assume in (4.11) that $n_1, n_2 > 0$. \\

This, however, contradicts (4.13), (4.14). \\  

{\bf Corollary 4.1} \\

A nontrivial linearly ordered monoid $( E, +, \leq )$ which is not a hyperspace, is Archimedean in the sense of (4.2). \\

\end{document}